\documentstyle[12pt]{article}

\def\bbbr{{\rm I\!R}}     %Real numbers
\def\bbbc{{\rm I\!\!\! C}}  % Complex numbers
\def\Fun{\rm Fun\, }
\def\ot{\otimes }

\def\Sym{\rm Sym\,}
\def\Vect{\rm Vect\,}
\def\bbbk{{\bf K}}
\def\uqs{ U_q(sl(2))}

\def\us{ U(sl(2))}

\def\ot2{\otimes 2}

\def\op{\oplus}

\newtheorem{proposition}{Proposition}

\newtheorem{theorem}{Theorem}
\newtheorem{definition}{Definition}
\newtheorem{remark}{Remark}
\def\adots{\mathinner{\mkern2mu\raise1pt\hbox{.}
\mkern3mu\raise4pt\hbox{.}\mkern1mu\raise7pt\hbox{.}}}

\title{Quantum anchor : $\uqs$ case}

\author{P.~Akueson\\
ISTV, Universit\'e de Valenciennes
59304 Valenciennes, France}

\begin{document}
\maketitle 
\begin{abstract}
We introduce the tangent space $\,T({\rm H}_q)\,$ on the quantum 
hyperboloid ($\,{\cal A}_{0,q}^c\,$) and equip it with an action  
 on $\,{\cal A}_{0,q}^c\,$ being a deformation of the action 
of vectors fields on functions. An embedding 
$\,\displaystyle\, sl(2)_q\,\longrightarrow\,T({\rm H}_q)\,$ of $q$-deformed Lie 
algebra $\,sl(2)\,$ being an analogue of the anchor 
$\,\displaystyle\, sl(2)\,\longrightarrow\,\Vect({\rm H})\,$ is called 
\lq\lq{quantum anchor}". 
\end{abstract}

\section{Introduction}

Let us consider the sphere of radius $\,R\,>\,0$
$$S^2\,=\,\{(x\,,\,y\,,\,z)\,\in\,{\bbbk}^3; \quad 
x^2 + y^2 + z^2\,=\,R^2\,\}$$
and the infinitesimal rotations 
$$X\,=\,y \partial_z - z\partial_y ,\quad Y\,=\,z \partial_x - x\partial_z,\quad
Z\,=\,x \partial_y - y\partial_x.$$
(Hereafter the basic field $\,{\bbbk}\,$ is $\,{\bbbr}\,$, but all results are 
still valid for its complexification. In this case we assume
$\,{\bbbk}={\bbbc}$.)

It is known that these infinitesimal rotations generate the space 
$\,\Vect(S^2)\,$ (or 
$\,\displaystyle\,\Vect({\rm A})\,,$ where 
$\,\displaystyle \,{\rm A}=\Fun(S^2)\,$ is the space of functions on the 
sphere) and  $\,\displaystyle\,\Vect(S^2)\,$ is the space of all vectors 
fields on $\,S^2$. (In what follows all functions are assumed to be 
restrictions of polynomials onto the sphere in question.) This means that every
element $\,{\cal X}\,\in\,\Vect(S^2)\,$ is of the form 
$${\cal X}\,=\,\alpha\,X + \beta\,Y + \gamma\,Z\,;\quad \alpha,\,\beta,\,
\gamma\,\in\,\Fun(S^2).$$
However, the vectors fields $\,X,\,Y,\,Z\,$ are not free. They satisfy the 
following identity
$$xX + yY + zZ\,=\,0.$$
So, the tangent space $\,T(S^2)\,$ on the sphere 
considered as a (left) $\,{\rm A}$-module can 
be defined by this identity, i.e., it is the factor-module of the rank $3$ free 
module 
$${\rm M}\,=\,\{\,\alpha\,X + \beta\,Y + \gamma\,Z\,;\quad  \alpha,\,\beta,\,
\gamma\,\in\,\Fun(S^2)\,\}$$
over its submodule 
$${\rm M}_1\,=\,\{\,f(xX + yY + zZ)\,;\quad f\,\in\,\Fun(S^2)\,\}.$$

Let us fix $\, k\,\in\,{\bbbk}\,$ and consider a hyperboloid
$${\rm H}\,=\,\{(u\,,\,v\,,\,w)\,\in\,{\bbbk}^3; \quad 2uw + vv + 2wu \,
=\,k\,\}.$$
Similarly to the sphere $\,S^2\,$, the  space of vectors fields on the hyperboloid  
$\,\Vect({\rm H})\,$ (or $\,\Vect({\rm A})\,$ where $\, A=\Fun({\rm H})\,$) is 
generated by three infinitesimal hyperbolic rotations $\,U,\,V,\,W\,$  
satisfying the identity
\begin{equation}
2uW + vV + 2wU\,=\,0. \label{eq1}
\end{equation}
Then, similarly to the previous case we can introduce the tangent space  
$\,T({\rm H})\,$ on the hyperboloid by means of the identity (\ref{eq1}).\\
We remark that in both cases the tangent modules are equipped with a Lie 
algebra structure and that an action
$$\beta \,:\,T({\rm H}) \otimes A\,\longrightarrow \,A\,,\qquad A\,=\,
\Fun({\rm H})$$
is well defined :  apply a vector field $\,{\cal X}\,$ to a 
function $\,f\,$ :
$${\cal X}\otimes f\,\longrightarrow \,{\cal X}f.$$
Thus, the tangent space $\,T({\rm H})\,$ is realized as that of vectors fields 
$\,\Vect({\rm H})\,$  on $\,{\rm H}\,$ 
and, consequently, we have an embedding of the Lie algebra 
$\,sl(2)\,$ generated by the elements $\,u,\,v,\,w\,$ into the space 
$\,\Vect({\rm H})\,$ generated by the elements $\,U,\,V,\,W\,$. 
Traditionally  such an embedding is called an {\it anchor}.\\

Consider  a quantum hyperboloid. The explicit description  
below, let us already give  three  proprieties  
of the corresponding \lq\lq{quantum function}" algebra : 
\begin{enumerate}
\item this algebra  is a flat deformation of its classical counterpart
(see  for example \cite{DGK} for the definition of a flat deformation), 
\item its product is $\,\uqs$-covariant (that is, it is a $\,\uqs$-morphism),
\item it is in some sens $q$-commutative (more precisely, it is a 
\lq\lq{$q$-commutative}" algebra of a two parameters family of $\,\uqs$-covariant 
algebras considered below).
\end{enumerate}

In the present note we discuss two problems : what is the \lq\lq{tangent space}" 
on the quantum hyperboloid ? Do we have an analogue of the above mentioned anchor ?\\
First of all, we should find  quantum analogues of the vectors fields 
$\,U,\,V,\,W$.
Usually, one considers the generators $\,X,\,Y,\,H\,$ 
of the quantum group (QG) $\,\uqs\,$ (see Section  1) as such quantum 
analogues. However, contrary to the classical case, those generators do not  
satisfy 
any $q$-analogue of identity (\ref{eq1}). So,
if we introduce the tangent space on the quantum
hyperboloid as the familly of all linear combinations of these operators
(with coefficients in the $q$-analogue of
$\,\Fun({\rm H})\,$, called {\it quantum hyperboloid} in what follows),  
then we do not have any flat
deformation of the classical tangent space $\,T({\rm H})$.\\
Here we propose other candidates to play the role of 
$q$-analogues of the operators $\,U,\,V,\,W\,,$  and satisfying 
an identity that can be viewed  as a $q$-analogue of (\ref{eq1}). Thus, the tangent space
$\,T({\rm H}_q)\,$ on the 
quantum hyperboloid $\,{\rm H}_q\,$ defined as the set of all  
linear combinations of these quantum vectors fields modulo this identity
will be  a deformation with respect to its classical counterpart.\\
Once the tangent space is defined, the following natural problem arises :
does there exist an action 
$$T({\rm H}_q)\otimes {\rm A}\,\longrightarrow\,{\rm A}$$
which is a deformation of the classical anchor ? 
Here $\,A\,$ designs the quantum analogue of the algebra 
$\,\Fun({\rm H})$. We construct such an action and give an interpretation of the
elements of $\,T({\rm H}_q)\,$ as 
{\it braided vectors fields}. We will obtain an embedding of $\,sl(2)_q\,$ 
into $\,T({\rm H}_q)\,$, called {\it quantum anchor}.
Here $\,sl(2)_q\,$ is the \lq\lq{braided}" analogue of $\,sl(2)\,$ and 
$\,T({\rm H}_q)\,$ is generated by 
braided analogues of the vectors fields $\,U,\,V,\,W$. \\
The  organization of this note is the following.
In the Section 1, we recall the construction of a quantum hyperboloid.
In section 2 we define the tangent space on a quantum hyperboloid and 
equip it with  an  action  on all \lq\lq{function}" on the   
quantum hyperboloid. Finally, we defined a {\it quantum anchor}.\\
Throughout the note, the parameter $\,q \in {\bbbk}\,$ is assumed to be 
generic. 

\section{Quantum hyperboloid}

In order to introduce a quantum hyperboloid, let us consider the QG $\,\uqs$ 
defined, as usual, by the 
generators $\,X,\,Y,\,H\,$ subject to the well known relations (cf. \cite{CP}). 
Let us fix a coproduct in $\,\uqs\,$ and the corresponding antipode and consider a
spin $1$ $\,\uqs$-module $\,{\rm V} = {\bf V}^q$.\\
We  only need the fact that the fusion ring of finite dimensional 
$\,\uqs$-modules is exactly
the same as in the classical case (we consider only the finite dimensional 
$\,\uqs$-modules wich are deformations of the classical $\,sl(2)$-modules). 
Thus, if $\,{\rm V}_i\,$ is a spin $\,i\,$ $\,\uqs$-module, the classical 
formula
$${\rm V}_i \otimes {\rm V}_j\,=\,\op_{k=\vert i-j \vert}^{i+j}\,{\rm V}_k$$
is still valid although the Clebsch-Gordan coefficients (which depend on the choice
of a base) are $q$-deformed.\\
In particular, we have :
$${\rm V}^{\otimes 2} \,=\,{\rm V}_0 \op {\rm V}_1 \op {\rm V}_2.$$
We keep the notation $\,{\rm V}\,$ for the initial space and 
we let $\,{\rm V}_1\,$ be   
the component in $\,{\rm V}^{\otimes 2}\,$ isomorphic to $\,{\rm V}$. Let us fix
in  the spaces $\,{\rm V},\,{\rm V}_0,\,{\rm V}_1\,$ and $\,{\rm V}_2\,$ some 
highest weight (h.w.) elements $\,v,\,v_0,\,v_1\,$ and $\,v_2\,$ respectively 
and impose the relations (which are the most general relations compatible with 
the action of the QG $\,\uqs\,$) :
\begin{equation}
v_0 - c\,=\,0 \qquad v_1 - \hbar\,v\,=\,0;\quad 
c,\,\hbar\,\in {\bbbk}. \label{eq2}
\end{equation}
One can now deduce the complete system of equations by  applying to these
relations the decreasing operator $\,Y\,\in\,\uqs$.\\ 
Let us denote $\,{\cal A}_{\hbar,q}^c\,$ the algebra defined by
the quotient algebra 
$$T({\rm V})\,/\,\{ I_{\hbar} \}$$
where $\,\displaystyle\,\{ I_{\hbar} \}\,$ is the ideal generated by the 
elements of the left hand side of  (\ref{eq2}) and all the derived elements.\\
The algebra $\,{\cal A}_{\hbar,q}^c\,$   is 
multiplicity free. More precisely, each integer spin module appears once in its 
decomposition into  a direct sum of irreducible $\,\uqs$-modules. Morever, each 
element of $\,{\cal A}_{\hbar,q}^c\,$ can be represented in a unique way as a 
sum of homogeneous elements such that each of them belongs to one of 
the highest components of 
$\,{\rm V}^{\otimes i}\,$ (a proof of this fact is given in \cite{A}).\\ 
We will treat the algebra  
$\,{\cal A}_{0,q}^c\,$ as a $q$-analogue of the commutative algebra 
$\,\Fun({\rm H})\,$ and call it  
{\it quantum hyperboloid} if  $\,c\not=0\,$ 
(being equipped with an involution it can  be treated as 
the  Podles quantum sphere (\cite{P})) and {\it quantum cone} if 
$\,c\,=\,0$. This is motivated by an analogy with the classical
case. Since the $\,sl(2)$-module $\,sl(2)^{\otimes 2}\,$ is
multiplicity free we can introduce $q$-analogues $\,I_{\pm}\,$ of
symmetric and skew symmetric subspaces of $\,sl(2)^{\otimes 2}\,$ by
setting, as in the classical case 
$$I_{+}\,=\,{\rm V}_0 \,\op \,{\rm V}_2 \qquad I_{-}\,=\,{\rm V}_1.$$
Let us emphasize that the corresponding algebras $\,{\rm A}_{\pm} =
{\rm T}(sl(2))/\{{I_{\mp}}\}\,$ are flat deformations of their classical
counterparts. Since the algebra $\,{\cal A}_{0,q}^c\,$ is a factor
of the \lq\lq{q-symmetric}" algebra $\,{\rm A}_{+}\,$ we  consider
it also as  \lq\lq{q-symmetric}". It is a particular case of the family of two 
parameters algebras $\,{\cal A}_{\hbar,q}^c\,$ (where $\,c\,$ is assumed to be fixed).

\noindent
\begin{remark}

In the case $\,n\geq 3\,$, the $\,sl(n)$-module $\,sl(n)^{\otimes
2}\,$ is not multiplicity free. Consequently, it is not obvious what  
$q$-analogues of the symmetric and skew symmetric  algebras of the
space $\,sl(n)\,$ should be. However, as shown in \cite{D}, there exists a flat deformation 
of the algebra $\,\Sym(sl(n))$. A way to construct these  
algebras in a more explicit way is suggested in \cite{AG}.\\
For any simple Lie algebra $\,g\,$ different from $\,sl(n)\,$
its tensor square $\,g^{\otimes 2}\,$ is multiplicity free. This implies that
there exists  natural $q$-analogues $\,I^q_{\pm}\,$ of the subspaces 
$\,I_{\pm} \subset g^{\otimes 2}$. However, the  $q$-algebras 
$\,{\rm A}^q_{\pm}=T({\rm V})/\{I^q_{\mp}\}\,$ are not flat deformations of 
their classical counterparts (cf. \cite{G}). 
\end{remark}

\section{Quantum anchor on quantum hyperboloid}

In the present section we introduce the tangent space on the quantum 
hyperboloid and  we equip it with a quantum anchor.\\
In the classical case, we can restate the identity (\ref{eq1}) in a symbolic 
way  as 
\begin{equation}
({\rm V}\otimes {\rm V}')_0 = 0 \label{eq4}
\end{equation}
where $\,{\rm V} = Span\,(u,\,v,\,w),\quad {\rm V}'=Span\,(U,\,V,\,W)\,$ and 
$\,\displaystyle\,({\rm V}\otimes {\rm V}')_i\,$ designs the 
spin $\,i\,$ $\,\us$-module  in the tensor product 
$\,\displaystyle\,{\rm V}\otimes {\rm V}'$. \\
In order to define the tangent space on the quantum hyperboloid, we should 
consider the identity (\ref{eq4})  in the $\,\uqs$-module category. Denote
$\,{\rm V}'^q\,$ the $q$-analogue of  $\,{\rm V}'$. Let  
$\,({\bf V}^q\otimes {\rm V}'^q)_0\,$ be the  spin $\,0\,$ $\,\uqs$-module  
 and $\,\{ ({\bf V}^q\otimes {\rm V}'^q)_0 \}\,$ be the left 
$\,{\cal A}_{0,q}^c$-module generated by  $\,({\bf V}^q \otimes {\rm V}'^q)_0$.

\begin{definition}
The \lq\lq{left tangent space}" on the quantum hyperboloid consider as a 
left $\,{\cal A}_{0,q}^c$-module is defined as follows :
\begin{equation}
T({\rm H}_q)_l = ({\cal A}_{0,q}^c \otimes {\rm V}'^q) /\{ ({\bf V}^q \otimes 
{\rm V}'^q)_0 \}. \label{eq5}
\end{equation}
\end{definition}
In a similar way one can define the  
right tangent space  $\,T({\rm H}_q)_r\,$  of the quantum hyperboloid 
considered as a right $\,{\cal A}_{0,q}^c$-module.

\begin{proposition}(\cite{A}, \cite{AG})

The $\,{\cal A}_{0,q}^c$-module  $\,T({\rm H}_q)\,$ is a flat deformation of
its classsical counterpart.
\end{proposition}

\noindent
Now we discuss the problem of a suitable definition of a \lq\lq{quantum anchor}"
on the quantum hyperboloid.\\
For this, we need the notion of a braided Lie algebra. Following  
\cite{DG}, we define a braided Lie bracket 
\begin{eqnarray*}
\displaystyle\,[ ,]_q\,:\,{\rm V}^{\otimes 2}\,\longrightarrow \,{\rm V}
\end{eqnarray*}
as a no trivial $\,\uqs$-morphism (here $\,{\rm V}\,$ is spin $1$ 
$\,\uqs$-module). By this request, the bracket is defined uniquely up to a 
factor. The product table of such a bracket is given for example in \cite{DG}.
Let us note that this construction is generalized for the $\,sl(n)\,$  $\,n>2\,$  case, 
in \cite{LS}.\\
The space $\,{\rm V}\,$ equipped with such bracket will be denoted by :
$$sl(2)_q\,:=\,({\rm V}\,,\,[\,,\,]_q)$$
and will be called a braided Lie algebra. The enveloping algebra 
of this braided Lie algebra can be defined similarly to that of 
$\,{\cal A}_{\hbar,q}^c\,$ but with ometting the first relation of  (\ref{eq2}). 
This enveloping algebra will be denoted $\,{\cal A}_{\hbar,q}$. We disregard here 
the problem of suitable relation between $\,\hbar\,$ and the factor coming 
in the definition of $\,sl(2)_q\,$, cf. \cite{LS}.\\
Let us consider (left) $q$-adjoint operators associated to the elements of the 
space $\,{\rm V}$. For example the operator $\,U^q\,$ asociated to 
the generator $\,u\,$ is defined by
$$  
\begin{array}{cccl} U^q :&{\rm V}&\longrightarrow&{\rm V}\\ 
&z&\mapsto&\,\displaystyle\,U^q z\,=\,ad^q u(z)\,=\,[u,z]_q ,\ \ U^q 1=0.
\end{array}$$
In a similar way the operators $\,V^q,\,W^q\,$  associated   
to the generators $\,v,\,w\,$  respectively are well defined on the space $\,{\rm V}$.

\begin{proposition}\label{pap}
The following holds
\begin{equation}
(q^3 +q )\,u\,W^q + v\,V^q + (q + q^{-1})\,w\,U^q\,=\,0.\label{eq6}
\end{equation}
Here all components are treated as operators acting from the space $\,{\rm V}\,$ to 
$\,{\cal A}_{0,q}^c$.
\end{proposition}

Throughout, we can assumed that the operators $\,U^q,\,V^q,\,W^q\,$ realize a 
representation of the braided Lie algebra $\,sl(2)_q\,$ (i.e. they satisfy 
the defining relation of its enveloping algebra). If it is not the case we can 
get such operators  by a proper rescaling 
$$U^q\,\longrightarrow\,\lambda\,U^q,\quad V^q\,\longrightarrow\,\lambda\,
V^q,\quad W^q\,\longrightarrow\,\lambda\,W^q \quad \lambda\,\in\,{\bbbk}.  $$
Let us remark that such a rescaling does not breack the identity (\ref{eq6}).

\begin{definition}\label{zuc}
We say that the $\,{\cal A}_{0,q}^c$-module $\,T({\rm H}_q)\,$ is equipped 
with a structure of \lq\lq{left quantum anchor}"  if there 
exists an action
$$\beta \,:\,T({\rm H}_q) \otimes {\cal A}_{0,q}^c\,\longrightarrow \,
{\cal A}_{0,q}^c\,$$
such taht the operators corresponding to $\,U^q,\,V^q,\,W^q\,$ realize a 
representation of the braided Lie algebra $\,sl(2)_q\,$ and the diagram  
$$
\begin{array}{ccc} 
{\cal A}_{0,q}^c \otimes T({\rm H}_q) \otimes {\cal A}_{0,q}^c&\longrightarrow 
 &T({\rm H}_q)\otimes {\cal A}_{0,q}^c\\  
\downarrow& &\downarrow\\ 
{\cal A}_{0,q}^c \otimes {\cal A}_{0,q}^c&\longrightarrow  &{\cal A}_{0,q}^c 
\end{array}
$$  
is commutative. Here we suppose that the elements of $\,{\cal A}_{0,q}^c\,$
act on $\,{\cal A}_{0,q}^c\,$ (the low arrow) by the usual product. The
vertical arrows are defined by  means of $\,\beta\,$ and the top one makes 
use of the $\,{\cal A}_{0,q}^c$-module structure of $\,T({\rm H}_q)$.
\end{definition}

\noindent
In a similar way a notion of \lq\lq{right quantum anchor}" can be 
defined.\\
Now we want to describe a quantum anchor in terms of the operators
$\,U^q,\,V^q,\,W^q$. However,
up to now the action of the operators $\,U^q,\,V^q,\,W^q\,$ is well 
defined on the space $\,{\rm V}$. To complet the construction of quantum anchor 
we should extend their action on the whole algebra $\,{\cal A}_{0,q}^c\,$ 
with proporties from definition \ref{zuc}. In the classical case, such an 
extension can be done via the Leibniz rule.\\
In the case under consideration such an extension is realize in \cite{A} 
(without any Leibniz rule). We represent here only the final result.

\begin{theorem}(\cite{A}) 
A quantum anchor exists and, moreover, it is unique if we impose the 
additional condition that
the extented operators $\,U^q,\,V^q,\,W^q\,$ send the highest component
($\,{\rm V}_i\,$) of $\,{\rm V}^{\otimes i}\,$ into itself. 
\end{theorem}

\noindent
To finish  this note, let us remark that the tangent space $\,T({\rm H}_q)\,$
is realize as an operator algebra on $\,{\cal A}_{0,q}^c\,$ and, moreover, the 
braided Lie algebra $\,sl(2)_q\,$ is embeded in it. This is the motivation
of our definition of  
quantum anchor, inspite of the fact that we do not equip $\,T({\rm H}_q)\,$ with 
any structure fo Lie algebra (deformed or \lq\lq{braided}").

 \vskip 1truecm

{\bf Acknowledgment :} I  would like to thank Prof  D. Gurevich for putting the
problem and  numerous helpful discussions.

\end{document}